\input amstex
\input amsppt.sty
 \magnification=\magstep1
 \hsize=30truecc
 \vsize=22.2truecm
 \baselineskip=16truept
 \nologo
 \TagsOnRight
 \def\N{\Bbb N}
 \def\Z{\Bbb Z}

 \def\C{\Bbb C}
 \def\l{\left}
 \def\r{\right}
 \def\bg{\bigg}
 \def\({\bg(}
 \def\[{\bg[}
 \def\){\bg)}
 \def\]{\bg]}
 \def\t{\text}
 \def\f{\frac}
 \def\colon{{:}\;}
 \def\ch{\roman{ch}}
 
 \def\per{\roman{per}}
 \def\Tor{\roman{Tor}}
 
 \def\se {\subseteq}
 
 \def\sm{\setminus}

 \def\bi{\binom}
 
 \def\cs{\ldots}
 \def\ls{\leqslant}
 \def\gs{\geqslant}
 
 \def\ve{\varepsilon}
 \def\da{\delta}
 
 \def\la{\lambda}

 \def\bi{\binom}
 
 \def\Proof{\noindent{\it Proof}}
 
 \def\Remark{\medskip\noindent{\it  Remark}}
 
\hbox{Math. Res. Lett. 15(2008), no.\,6, 1263-1276}
\bigskip
  \topmatter
  \title An additive theorem and restricted sumsets\endtitle
  \author Zhi-Wei Sun \endauthor
 \address Department of Mathematics, Nanjing University,
   Nanjing 210093, People's Republic of China\endaddress
 \email {zwsun\@nju.edu.cn}\ {\it Homepage}:\ {\tt http://math.nju.edu.cn/$\sim$zwsun}
 \endemail

 \abstract
 Let $G$ be any additive abelian group with cyclic torsion subgroup,
and let $A$, $B$ and $C$ be finite subsets of $G$ with cardinality
$n>0$. We show that there is a numbering $\{a_i\}_{i=1}^n$ of the
elements of $A$, a numbering $\{b_i\}_{i=1}^n$ of the elements of
$B$ and a numbering $\{c_i\}_{i=1}^n$ of the elements of $C$, such
that all the sums $a_i+b_i+c_i\ (1\ls i\ls n)$ are (pairwise)
distinct. Consequently, each subcube of the Latin cube formed by
the Cayley addition table of $\Z/N\Z$ contains a Latin
transversal. This additive theorem is an essential result which
can be further extended via restricted sumsets in a field.
 \endabstract
 \thanks 2000 {\it Mathematics Subject Classification.}
Primary 11B75; Secondary 05A05, 05B15, 05E99, 11C08, 11P99, 15A15,
20D60, 20K99.
\newline
\indent Supported by the National Science Fund (grant 10425103) for
Distinguished Young Scholars in China.
 \endthanks
 \endtopmatter
  \document
 \heading {1. Introduction}\endheading

 In 1999 Snevily [Sn] raised the following beautiful conjecture in additive combinatorics
 which is currently an active area of research.

\proclaim{Snevily's Conjecture} Let $G$ be an additive abelian
group with $|G|$ odd. Let $A$ and $B$ be subsets of $G$ with
cardinality $n\in\Z^+=\{1,2,3,\ldots\}$. Then there is a numbering
$\{a_i\}_{i=1}^n$ of the elements of $A$ and a numbering
$\{b_i\}_{i=1}^n$ of the elements of $B$ such that the sums
$a_1+b_1,\ldots,a_n+b_n$ are $($pairwise$)$ distinct.
\endproclaim

When $|G|$ is an odd prime, this conjecture was proved by Alon
[A2] via the polynomial method rooted in Alon and Tarsi [AT], and
developed by Alon, Nathanson and Ruzsa [ANR] (see also [N, pp.\,98-107]
and [TV, pp.\,329-345])
and refined by Alon [A1] in 1999. In 2001 Dasgupta,
K\'arolyi, Serra and Szegedy [DKSS] confirmed Snevily's conjecture
for any cyclic group of odd order. In 2003 Sun [Su3] obtained some
further extensions of the Dasgupta-K\'arolyi-Serra-Szegedy result via restricted sums in a field.

In Snevily's conjecture the abelian group is required to have odd
order. (An abelian group of even order has an element $g$ of order
2 and hence we don't have the described result for $A=B=\{0,g\}$.)
For a general abelian group $G$ with its torsion subgroup $\Tor(G)=\{a\in G:\, a\ \t{has a finite order}\}$ cyclic, if
we make no hypothesis on the order of $G$, what additive
properties can we impose on several finite subsets of $G$ with
cardinality $n$? In this direction we establish the following new
theorem of additive nature.

\proclaim{Theorem 1.1} Let $G$ be any additive abelian group with cyclic torsion subgroup,
and let $A_1,\ldots,A_m$
be arbitrary subsets of $G$ with cardinality $n\in\Z^+$, where $m$ is odd.
Then the elements of $A_i\ (1\ls i\ls m)$ can be listed in a suitable order $a_{i1},\ldots,a_{in}$,
so that all the sums $\sum_{i=1}^ma_{ij}\ (1\ls j\ls n)$ are distinct.
In other words, for a certain subset $A_{m+1}$ of $G$ with $|A_{m+1}|=n$,
there is a matrix $(a_{ij})_{1\ls i\ls m+1,\,1\ls j\ls n}$ such that $\{a_{i1},\ldots,a_{in}\}=A_i$
for all $i=1,\ldots,m+1$ and the column sum $\sum_{i=1}^{m+1}a_{ij}$ vanishes for every $j=1,\ldots,n$.
\endproclaim

\Remark\ 1.1. Theorem 1.1 in the case $m=3$ is essential;
the result for $m=5,7,\ldots$ can be obtained by repeated use of the case $m=3$.

\medskip

\noindent{\bf Example 1.1}. In Theorem 1.1 the condition $2\nmid m$ is indispensable.
Let $G$ be an additive cyclic group of even order $n$. Then $G$ has a unique element $g$ of order 2
and hence $a\not=-a$ for all $a\in G\sm\{0,g\}$. Thus $\sum_{a\in G}a=0+g=g$.
For each $i=1,\ldots,m$ let $a_{i1},\ldots,a_{in}$ be a list of the $n$ elements of $G$.
If those $\sum_{i=1}^ma_{ij}$ with $1\ls j\ls n$ are distinct, then
$$\sum_{a\in G}a=\sum_{j=1}^n\sum_{i=1}^ma_{ij}=\sum_{i=1}^m\sum_{j=1}^na_{ij}=m\sum_{a\in G}a,$$
hence $(m-1)g=(m-1)\sum_{a\in G}a=0$ and therefore $m$ is odd.

\medskip

\noindent{\bf Example 1.2}. The group $G$ in Theorem 1.1 cannot be replaced
by an arbitrary abelian group. To illustrate this, we look at the Klein quaternion group
$$\Z/2\Z\oplus\Z/2\Z=\{(0,0),(0,1),(1,0),(1,1)\}$$
and its subsets
$$A_1=\{(0,0),(0,1)\},\, A_2=\{(0,0),(1,0)\},\, A_3=\cdots=A_m=\{(0,0),(1,1)\},$$
where $m\gs3$ is odd. For $i=1,\ldots,m$ let $a_i,a_i'$ be a list of the two elements of $A_i$,
then
$$\sum_{i=1}^m(a_i+a_i')=(0,1)+(1,0)+(m-2)(1,1)=(0,0)$$
and hence
$\sum_{i=1}^ma_i=-\sum_{i=1}^ma_i'=\sum_{i=1}^ma_i'$.
\medskip

Recall that a line of an $n\times n$ matrix is a row or column of the matrix.
We define a line of an $n\times n\times n$ cube in a similar way.
A {\it Latin cube} over a set $S$ of cardinality $n$ is an $n\times n\times n$ cube
whose entries come from the set $S$ and
no line of which contains a repeated element. A {\it transversal} of an $n\times n\times n$ cube
is a collection of $n$ cells no two of which lie in the same line. A {\it Latin transversal}
of a cube is a transversal whose cells contain no repeated element.

\proclaim{Corollary 1.1} Let $N$ be any positive integer. For the
$N\times N\times N$ Latin cube over $\Z/N\Z$ formed by the Cayley addition table,
each $n\times n\times n$ subcube with $n\ls N$ contains a Latin transversal.
\endproclaim
\Proof. Just apply Theorem 1.1 with $G=\Z/N\Z$ and $m=3$. \qed

\smallskip
In 1967 Ryser [R] conjectured that every Latin square of odd order
has a Latin transversal. Another conjecture of Brualdi (cf. [D],
[DK, p.\,103] and [EHNS]) states that every Latin square of order
$n$ has a partial Latin transversal of size $n-1$. These and
Corollary 1.1 suggest that our following conjecture might be
reasonable.

\proclaim{Conjecture 1.1} Every $n\times n\times n$ Latin cube contains a Latin transversal.
\endproclaim

Note that Conjecture 1.1 does not imply Theorem 1.1 since an $n\times n\times n$ subcube of a Latin cube
might have more than $n$ distinct entries.

\proclaim{Corollary 1.2} Let $G$ be any additive abelian group with cyclic torsion subgroup,
and let $A_1,\ldots,A_m$
be subsets of $G$ with cardinality $n\in\Z^+$, where $m$ is even.
Suppose that all the elements of $A_m$ have odd order.
Then the elements of $A_i\ (1\ls i\ls m)$ can be listed in a suitable order $a_{i1},\ldots,a_{in}$,
so that all the sums $\sum_{i=1}^ma_{ij}\ (1\ls j\ls n)$ are distinct.
\endproclaim
\Proof. As $m-1$ is odd, by Theorem 1.1 the elements of $A_i\
(1\ls i\ls m-1)$ can be listed in a suitable order
$a_{i1},\ldots,a_{in}$, such that all the sums
$s_j=\sum_{i=1}^{m-1}a_{ij}\ (1\ls j\ls n)$ are distinct. Since
all the elements of $A_m$ have odd order, by [Su3, Theorem
1.1(ii)] there is a numbering $\{a_{mj}\}_{j=1}^n$ of the elements
of $A_m$ such that all the sums $s_j+a_{mj}=\sum_{i=1}^ma_{ij}\
(1\ls j\ls n)$ are distinct. We are done. \qed

\smallskip

As an essential result, Theorem 1.1 might have various potential applications
in additive number theory and combinatorial designs.

We can extend Theorem 1.1 via restricted sumsets in a field.
The additive order of the multiplicative
identity of a field $F$ is either infinite or a prime;
we call it the {\it characteristic} of $F$ and denote it by $\ch(F)$.
The reader is referred to [DH], [ANR], [Su2], [HS],
[LS], [PS1], [Su3], [SY] and [PS2] for various results on restricted sumsets of the type
$$\{a_1+\cdots+a_n:\ a_1\in A_1,\ldots,a_n\in A_n\ \t{and}\ P(a_1,\ldots,a_n)\not=0\},$$
where $A_1,\ldots,A_n\se F$ and $P(x_1,\ldots,x_n)\in F[x_1,\ldots,x_n]$.

For a finite sequence $\{A_i\}_{i=1}^n$ of sets, if $a_1\in
A_1,\ldots,a_n\in A_n$ and $a_1,\ldots,a_n$ are distinct, then the
sequence $\{a_i\}_{i=1}^n$ is called a {\it system of distinct
representives} (SDR) of $\{A_i\}_{i=1}^n$. This concept plays an
important role in combinatorics and a celebrated theorem of Hall
tells us when $\{A_i\}_{i=1}^n$ has an SDR (see, e.g., [Su1]).
Most results in our paper involve SDRs of several subsets of a
field.

Now we state our second theorem which is much more general than Theorem 1.1.

\proclaim{Theorem 1.2} Let $h,k,l,m,n$ be positive integers satisfying
$$k-1\gs m(n-1)\ \ \t{and}\ \ l-1\gs h(n-1).\tag1.1$$
Let $F$ be a field with $\ch(F)>\max\{K,L\}$, where
$$K=(k-1)n-(m+1)\bi n2\ \t{and}\ L=(l-1)n-(h+1)\bi n2.\tag1.2$$
Assume that $c_1,\ldots,c_n\in F$ are distinct and
$A_1,\ldots,A_n,B_1,\ldots,B_n$ are subsets of $F$ with
$$|A_1|=\cdots=|A_n|=k\ \t{and}\ |B_1|=\cdots=|B_n|=l.\tag1.3$$
Let $P_1(x),\ldots,P_n(x),Q_1(x),\ldots,Q_n(x)\in F[x]$ be monic polynomials
with $\deg P_i(x)=m$ and $\deg Q_i(x)=h$ for $i=1,\ldots,n$.
Then, for any $S,T\se F$ with $|S|\ls K$ and
$|T|\ls L$, there exist $a_1\in A_1,\ldots,a_n\in A_n,b_1\in B_1,\ldots,b_n\in B_n$ such that
$a_1+\cdots+a_n\not\in S$, $b_1+\cdots+b_n\not\in T$, and also
$$a_ib_ic_i\not=a_jb_jc_j,\ P_i(a_i)\not=P_j(a_j),\ Q_i(b_i)\not=Q_j(b_j)\ \ \t{if}\ 1\ls i<j\ls n.\tag1.4$$
\endproclaim

\Remark\ 1.2. If $h,k,l,m,n$ are positive integers satisfying (1.1), then the integers $K$ and $L$ given by (1.2)
are nonnegative since
$$K\gs m(n-1)n-(m+1)\bi n2=(m-1)\bi n2\ \t{and}\ L\gs (h-1)\bi n2.$$
\medskip

 From Theorem 1.2 we can deduce the following extension of Theorem 1.1.

 \proclaim{Theorem 1.3} Let $G$ be an additive abelian group with cyclic torsion subgroup.
 Let $h,k,l,m,n$ be positive integers satisfying $(1.1)$.
 Assume that $c_1,\ldots,c_n\in G$ are distinct, and $A_1,\cs,A_n,B_1,\ldots,B_n$ are
 subsets of $G$ with $|A_1|=\cdots=|A_n|=k$ and $|B_1|=\cdots=|B_n|=l$.
 Then, for any sets $S$ and $T$ with $|S|\ls (k-1)n-(m+1)\bi n2$ and $|T|\ls (l-1)n-(h+1)\bi n2$, there are
 $a_1\in A_1,\ldots,a_n\in A_n,b_1\in B_1,\ldots,b_n\in B_n$ such that
 $\{a_1,\ldots,a_n\}\not\in S$, $\{b_1,\ldots,b_n\}\not\in T$, and also
$$a_i+b_i+c_i\not=a_j+b_j+c_j,\ ma_i\not=ma_j,\ hb_i\not=hb_j\ \ \t{if}\ 1\ls i<j\ls n.\tag1.5$$
\endproclaim
 \Proof. Let $H$ be the subgroup of $G$ generated by the finite set
 $$A_1\cup\cdots\cup A_n\cup B_1\cup\cdots\cup B_n\cup\{c_1,\ldots,c_n\}.$$
 Since $\Tor(H)$ is cyclic and finite, as in the proof of [Su3, Theorem 1.1]
 we can identify the additive group $H$ with a subgroup of the multiplicative group
 $\C^*=\C\sm\{0\}$, where $\C$ is the field of complex numbers.
 So, without loss of generality, below we simply view $G$ as the
 multiplicative group $\C^*$.

 Let $S$ and $T$ be two sets with  $|S|\ls (k-1)n-(m+1)\bi n2$ and $|T|\ls (l-1)n-(h+1)\bi
 n2$. Then
 $$S'=\{a_1+\cdots+a_n:\ a_1\in A_1,\ldots,a_n\in A_n,\
 \{a_1,\ldots,a_n\}\in S\}$$
 and
$$T'=\{b_1+\cdots+b_n:\ b_1\in B_1,\ldots,b_n\in B_n,\
 \{b_1,\ldots,b_n\}\in T\}$$
 are subsets of $\C$ with $|S'|\ls |S|$ and $|T'|\ls |T|$.
 By Theorem 1.2 with $P_i(x)=x^m$ and $Q_i(x)=x^h\ (1\ls i\ls
 n)$, there are
 $a_1\in A_1,\ldots,a_n\in A_n,b_1\in B_1,\ldots,b_n\in B_n$ such that
 $a_1+\cdots+a_n\not\in S'$ (and hence $\{a_1,\ldots,a_n\}\not\in S$), $b_1+\cdots+b_n\not\in T'$
 (and hence $\{b_1,\ldots,b_n\}\not\in T$), and also
 $$a_ib_ic_i\not=a_jb_jc_j, \ a_i^m\not=a_j^m,\ b_i^h\not=b_j^h\ \ if\ 1\ls i<j\ls n.$$
This concludes the proof.  \qed

 \Remark\ 1.3. Theorem 1.1 in the case $m=3$ is a special case of Theorem 1.3.
 \medskip

 Here is another extension of Theorem 1.1 via restricted sumsets in a field.

 \proclaim{Theorem 1.4} Let $k,m,n$ be positive integers with $k-1\gs m(n-1)$, and let
 $F$ be a field with $\ch(F)>\max\{mn,(k-1-m(n-1))n\}$. Assume that $c_1,\ldots,c_n\in F$ are
 distinct, and $A_1,\ldots,A_n,B_1,\ldots,B_n$ are subsets of $F$ with $|A_1|=\cdots=|A_n|=k$
 and $|B_1|=\cdots=|B_n|=n$. Let $S_{ij}\se F$ with $|S_{ij}|<2m$ for all $1\ls i<j\ls n$.
 Then there is an SDR $\{b_i\}_{i=1}^n$ of $\{B_i\}_{i=1}^n$ such that the restricted sumset
 $$S=\{a_1+\cdots+a_n:\, a_i\in A_i,\ a_i-a_j\not\in S_{ij}\ \t{and}\ a_ib_ic_i\not=a_jb_jc_j\ \t{if}\ i<j\}\tag1.6$$
 has at least $(k-1-m(n-1))n+1$ elements.
 \endproclaim

 Now we introduce some basic notations in this paper.
Let $R$ be any commutative ring with
identity. The {\it permanent} of a matrix
$A=(a_{ij})_{1\ls i,j\ls n}$ over $R$ is given by
$$\per(A)=\|a_{ij}\|_{1\ls i,j\ls n}=\sum_{\sigma\in S_n}a_{1,\sigma(1)}\cdots a_{n,\sigma(n)},\tag1.7$$
where $S_n$ is the symmetric group of all the
permutations on $\{1,\cs,n\}$.
Recall that the determinant
of $A$ is defined by
$$\det(A)=|a_{ij}|_{1\ls i,j\ls n}=\sum_{\sigma\in S_n}\ve(\sigma)a_{1,\sigma(1)}\cdots a_{n,\sigma(n)},\tag1.8$$
where $\ve(\sigma)$ is $1$ or $-1$ according as $\sigma$ is even
or odd. We remind the difference between the notations $|\cdot|$
and $\|\cdot\|$. For the sake of convenience, the coefficient of
the monomial $x_1^{k_1}\cdots x_n^{k_n}$ in a polynomial
$P(x_1,\ldots,x_n)$ over $R$ will be denoted by $[x_1^{k_1}\cdots
x_n^{k_n}]P(x_1,\ldots,x_n)$.

 In the next section we are going to prove Theorem 1.1 in two different ways.
  Section 3 is devoted to the study of duality between determinant and permanent.
 On the basis of Section 3, we will show Theorem 1.2 in Section 4 via the polynomial method.
In Section 5, we will present our proof of Theorem 1.4.

\heading{2. Two proofs of Theorem 1.1}\endheading

\proclaim{Lemma 2.1} Let $R$ be a commutative ring with identity,
and let $a_{ij}\in R$ for $i=1,\ldots,m$ and $j=1,\ldots,n$,
where $m\in\{3,5,\ldots\}$. The we have the identity
$$\aligned&\sum_{\sigma_1,\ldots,\sigma_{m-1}\in S_n}\ve(\sigma_1\cdots\sigma_{m-1})
\prod_{1\ls i<j\ls
n}\(a_{mj}\prod_{s=1}^{m-1}a_{s\sigma_s(j)}-a_{mi}\prod_{s=1}^{m-1}a_{s\sigma_s(i)}\)
\\&\qquad\qquad=\prod_{1\ls i<j\ls n}(a_{1j}-a_{1i})\cdots(a_{mj}-a_{mi}).
\endaligned\tag2.1$$
\endproclaim
\Proof. Recall that $|x_j^{i-1}|_{1\ls i,j\ls n}=\prod_{1\ls i<j\ls n}(x_j-x_i)$ (Vandermonde).
Let $\Sigma$ denote the left-hand side of (2.1). Then
$$\align\Sigma=&\sum_{\sigma_1,\ldots,\sigma_{m-1}\in S_n}\ve(\sigma_1\cdots\sigma_{m-1})
|(a_{1,\sigma_1(j)}\cdots a_{m-1,\sigma_{m-1}(j)}a_{mj})^{i-1}|_{1\ls i,j\ls n}
\\=&\sum_{\sigma_1,\ldots,\sigma_{m-1}\in S_n}\ve(\sigma_1)\times\cdots\times\ve(\sigma_{m-1})
\\&\quad\ \ \ \times\sum_{\tau\in S_n}\ve(\tau)\prod_{i=1}^n
(a_{1,\sigma_1(\tau(i))}\cdots a_{m-1,\sigma_{m-1}(\tau(i))}a_{m,\tau(i)})^{i-1}
\\=&\sum_{\tau\in S_n}\ve(\tau)^m\prod_{i=1}^na_{m,\tau(i)}^{i-1}
\times\prod_{s=1}^{m-1}\sum_{\sigma_s\in
S_n}\ve(\sigma_s\tau)\prod_{i=1}^na_{s,\sigma_s\tau(i)}^{i-1}
\\=&\sum_{\tau\in S_n}\ve(\tau)^m\prod_{i=1}^na_{m,\tau(i)}^{i-1}
\times\prod_{s=1}^{m-1}\sum_{\sigma\in
S_n}\ve(\sigma)\prod_{i=1}^na_{s,\sigma(i)}^{i-1}.
\endalign$$
Since $m$ is odd, we finally have
$$\align\Sigma=&|a_{mj}^{i-1}|_{1\ls i,j\ls n}\prod_{s=1}^{m-1}|a_{sj}^{i-1}|_{1\ls i,j\ls n}
=\prod_{s=1}^m\prod_{1\ls i<j\ls n}(a_{sj}-a_{si}).
\endalign$$
This proves (2.1). \qed

\Remark\ 2.1. When $m\in\{2,4,6,\ldots\}$, the right-hand side of (2.1) should be replaced by
$$\|a_{mj}^{i-1}\|_{1\ls i,j\ls n}\prod_{1\ls i<j\ls n}(a_{1j}-a_{1i})\cdots(a_{m-1,j}-a_{m-1,i}).$$
\medskip

\noindent{\bf Definition 2.1}. A subset $S$ of a commutative ring
$R$ with identity is said to be {\it regular} if all those $a-b$
with $a,b\in S$ and $a\not=b$ are units (i.e., invertible elements)
of $R$.

\proclaim{Theorem 2.1} Let $R$ be a commutative ring with
identity, and let $m>0$ be odd. Then, for any regular subsets
$A_1,\ldots,A_m$ of $R$ with cardinality $n\in\Z^+$, the elements
of $A_i\ (1\ls i\ls m)$ can be listed in a suitable order
$a_{i1},\ldots,a_{in}$, so that all the products
$\prod_{i=1}^ma_{ij}\ (1\ls j\ls n)$ are distinct.
\endproclaim
\Proof. The case $m=1$ is trivial. Below we let $m\in\{3,5,\ldots\}$.

Write $A_s=\{b_{s1},\ldots,b_{sn}\}$ for $s=1,\ldots,m$. As
all those $b_{sj}-b_{si}$ with $1\ls s\ls m$ and $1\ls i<j\ls n$ are units of $R$, the product
$$\prod_{1\ls i<j\ls n}(b_{1j}-b_{1i})\cdots(b_{mj}-b_{mi})$$
is also a unit of $R$ and hence nonzero.
Thus, by Lemma 2.1 there are $\sigma_1,\ldots,\sigma_{m-1}\in S_n$ such that whenever $1\ls i<j\ls n$ we have
$$b_{1,\sigma_1(i)}\cdots b_{m-1,\sigma_{m-1}(i)}b_{mi}\not=b_{1,\sigma_1(j)}\cdots b_{m-1,\sigma_{m-1}(j)}b_{mj}.$$
For $1\ls s\ls m$ and $1\ls j\ls n$, let
$a_{sj}=b_{s,\sigma_s(j)}$ if $s<m$, and $a_{sj}=b_{sj}$ if $s=m$.
Then $\{a_{s1},\ldots,a_{sn}\}=A_s$, and all the products
$\prod_{s=1}^ma_{sj}\ (j=1,\ldots,n)$ are distinct. This concludes
the proof. \qed

\medskip
\noindent{\it Proof of Theorem 1.1}. As mentioned in the proof of
Theorem 1.3 via Theorem 1.2, without loss of generality we may
simply take $G$ to be the multiplicative group $\C^*=\C\sm\{0\}$.
As any nonzero element of a field is a unit in the field, the
desired result follows from Theorem 2.1 immediately. \qed

\smallskip

Now we turn to our second approach to Theorem 1.1.

\proclaim{Lemma 2.2} Let $c_1,\ldots,c_n$ be elements of a commutative ring with identity.
Then we have
$$\aligned&[x_1^{n-1}\cdots x_n^{n-1}y_1^{n-1}\cdots y_n^{n-1}]
\prod_{1\ls i<j\ls n}(x_j-x_i)(y_j-y_i)(c_jx_jy_j-c_ix_iy_i)
\\&\qquad\qquad\qquad\qquad\qquad=\prod_{1\ls i<j\ls n}(c_j-c_i).\endaligned\tag2.2$$
\endproclaim
\Proof. Observe that
$$\align&\prod_{1\ls i<j\ls n}(x_j-x_i)(y_j-y_i)(c_jx_jy_j-c_ix_iy_i)
\\=&|x_i^{j-1}|_{1\ls i,j\ls n}|y_i^{j-1}|_{1\ls i,j\ls n}|(c_ix_iy_i)^{j-1}|_{1\ls i,j\ls n}
\\=&\sum_{\sigma\in S_n}\ve(\sigma)\prod_{i=1}^n x_i^{\sigma(i)-1}
\times\sum_{\tau\in S_n}\ve(\tau)\prod_{i=1}^n y_i^{\tau(i)-1}
\times\sum_{\la\in S_n}\ve(\la)\prod_{i=1}^n (c_ix_iy_i)^{\la(i)-1}
\\=&\sum_{\la\in S_n}\ve(\la)\prod_{i=1}^nc_i^{\la(i)-1}\sum_{\sigma,\tau\in S_n}\ve(\sigma\tau)\prod_{i=1}^n
\l(x_i^{\la(i)+\sigma(i)-2}y_i^{\la(i)+\tau(i)-2}\r).
\endalign$$
Thus the left-hand side of (2.2) coincides with
$$\sum_{\la\in S_n}\(\ve(\la)\prod_{i=1}^nc_i^{\la(i)-1}\)\ve(\bar\la\bar\la)=|c_i^{j-1}|_{1\ls i,j\ls n}
=\prod_{1\ls i<j\ls n}(c_j-c_i),$$ where $\bar\la(i)=n+1-\la(i)$
for $i=1,\ldots,n$. We are done. \qed
\medskip

Let us recall the following central principle of the polynomial method.

\proclaim{Combinatorial Nullstellensatz {\rm [A1]}} Let
$A_1,\ldots,A_n$ be finite subsets of a field $F$ with $|A_i|>k_i$
for $i=1,\ldots,n$, where $k_1,\ldots,k_n$ are nonnegative
integers. If the total degree of $f(x_1,\ldots,x_n)\in
F[x_1,\ldots,x_n]$ is $k_1+\cdots+k_n$ and $[x_1^{k_1}\cdots
x_n^{k_n}]f(x_1,\ldots,x_n)$ is nonzero, then
$f(a_1,\ldots,a_n)\not=0$ for some $a_1\in A_1,\ldots,a_n\in A_n$.
\endproclaim

\proclaim{Theorem 2.2} Let $A_1,\ldots,A_n$ and $B_1,\ldots,B_n$ be subsets of a field $F$ with cardinality $n$.
And let $c_1,\ldots,c_n$ be distinct elements of $F$. Then there is an SDR $\{a_i\}_{i=1}^n$ of
$\{A_i\}_{i=1}^n$ and an SDR $\{b_i\}_{i=1}^n$ of $\{B_i\}_{i=1}^n$
such that the products $a_1b_1c_1,\ldots,a_nb_nc_n$ are distinct.
\endproclaim
\Proof. As $c_1,\ldots,c_n$ are distinct, (2.2) implies that
$$[x_1^{n-1}\cdots x_n^{n-1}y_1^{n-1}\cdots y_n^{n-1}]
\prod_{1\ls i<j\ls n}(x_j-x_i)(y_j-y_i)(c_jx_jy_j-c_ix_iy_i)\not=0.$$
Applying the Combinatorial Nullstellensatz, we obtain the desired result. \qed

\Remark\ 2.2. When $F=\C$, $A_1=\cdots=A_n$ and $B_1=\cdots=B_n$,
Theorem 2.2 yields Theorem 1.1 with $m=3$. Note also that Theorems
1.2 and 1.4 are different extensions of Theorem 2.2.

\heading{3. Duality between determinant and permanent}\endheading

Let us first summarize Theorem 2.1 and Corollary 2.1 of Sun [Su3] in the following theorem.

\proclaim{Theorem 3.1 {\rm (Sun [Su3])}} Let $R$ be a commutative
ring with identity, and let $A=(a_{ij})_{1\ls i,j\ls n}$ be a
matrix over $R$.

{\rm (i)} Let $k_1,\ldots,k_n,m_1,\ldots,m_n\in\N=\{0,1,2,\ldots\}$ with
$M=\sum_{i=1}^n m_i+\da\bi n2\ls \sum_{i=1}^nk_i$ where $\da\in\{0,1\}$.
Then
$$\align&[x_1^{k_1}\cdots x_n^{k_n}]
|a_{ij}x_j^{m_i}|_{1\ls i,j\ls n}\prod_{1\ls i<j\ls n}
(x_j-x_i)^{\da}\times\bg(\sum_{s=1}^nx_s\bg)^{\sum_{i=1}^n k_i-M}
\\&\quad\ \ =\cases\sum_{\sigma\in S_n,\,
D_{\sigma}\se\N}\ve(\sigma)N_{\sigma}\prod_{i=1}^na_{i,\sigma(i)}&\t{if}\ \da=0,
\\\sum_{\sigma\in T_n}\ve(\sigma')N_{\sigma}
\prod_{i=1}^na_{i,\sigma(i)}&\t{if}\ \da=1,
\endcases\endalign$$
where
$$\align D_{\sigma}=&\{k_{\sigma(1)}-m_1,\cs,k_{\sigma(n)}-m_n\},
\\T_n=&\{\sigma\in S_n\colon D_{\sigma}\se\N\ \t{and}\ |D_{\sigma}|=n\},
\\N_{\sigma}=&\f{(k_1+\cdots+k_n-M)!}{\prod_{i=1}^n\prod\Sb0\ls j<k_{\sigma(i)}-m_i
\\ j\not\in D_{\sigma}\ \t{if}\ \da=1\endSb(k_{\sigma(i)}-m_i-j)}\in\Z^+,
\endalign$$
and $\sigma'\ ($with $\sigma\in T_n)$ is the unique permutation in $S_n$ such that
$$0\ls k_{\sigma(\sigma'(1))}-m_{\sigma'(1)}<\cdots<k_{\sigma(\sigma'(n))}-m_{\sigma'(n)}.$$

{\rm (ii)} Let $k,m_1,\cs,m_n\in\N$ with $m_1\ls\cdots\ls m_n\ls k$.
Then
$$\aligned&[x_1^k\cdots x_n^k]|a_{ij}x_j^{m_i}|_{1\ls i,j\ls
n}(x_1+\cdots+x_n)^{kn-\sum_{i=1}^n m_i}
\\&\qquad=\f{(kn-\sum_{i=1}^n m_i)!}{\prod_{i=1}^n(k-m_i)!}\det(A).
\endaligned\tag3.1$$
In the case $m_1<\cdots<m_n$, we also have
$$\aligned &[x_1^k\cdots x_n^k]|a_{ij}x_j^{m_i}|_{1\ls i,j\ls n}
\prod_{1\ls i<j\ls n}(x_j-x_i)\times\bg(\sum_{s=1}^nx_s\bg)^{kn-\bi n2-\sum_{i=1}^n m_i}
\\&\quad\ =(-1)^{\bi n2}\f{(kn-\bi n2-\sum_{i=1}^n m_i)!}
{\prod_{i=1}^n\prod\Sb m_i<j\ls k\\j\not\in\{m_s:\, i<s\ls
n\}\endSb(j-m_i)}\per(A).
\endaligned\tag3.2$$
\endproclaim

In view of the minor difference between the definitions of determinant and
permanent, by modifying the proof of the above result in [Su3] slightly
we get the following dual of Theorem 3.1.

\proclaim{Theorem 3.2} Let $R$ be a commutative ring with
identity, and let $A=(a_{ij})_{1\ls i,j\ls n}$ be a matrix over
$R$.

{\rm (i)} Let $k_1,m_1,\ldots,k_n,m_n\in\N$ with
$M=\sum_{i=1}^n m_i+\da\bi n2\ls \sum_{i=1}^nk_i$ where $\da\in\{0,1\}$.
Then
$$\align&[x_1^{k_1}\cdots x_n^{k_n}]
\|a_{ij}x_j^{m_i}\|_{1\ls i,j\ls n}\prod_{1\ls i<j\ls n}
(x_j-x_i)^{\da}\times\bg(\sum_{s=1}^nx_s\bg)^{\sum_{i=1}^n k_i-M}
\\&\quad\ \ =\cases\sum_{\sigma\in S_n,\,
D_{\sigma}\se\N}N_{\sigma}\prod_{i=1}^na_{i,\sigma(i)}&\t{if}\ \da=0,
\\\sum_{\sigma\in T_n}\ve(\sigma\sigma')N_{\sigma}
\prod_{i=1}^na_{i,\sigma(i)}&\t{if}\ \da=1,
\endcases\endalign$$
where $D_{\sigma},T_n, N_{\sigma}$ and $\sigma'$ are as in Theorem
{\rm 3.1(i)}.

{\rm (ii)} Let $k,m_1,\cs,m_n\in\N$ with $m_1\ls\cdots\ls m_n\ls k$.
Then
$$\aligned&[x_1^k\cdots x_n^k]\|a_{ij}x_j^{m_i}\|_{1\ls i,j\ls n}(x_1+\cdots+x_n)^{kn-\sum_{i=1}^n m_i}
\\&\qquad=\f{(kn-\sum_{i=1}^n m_i)!}{\prod_{i=1}^n(k-m_i)!}\per(A).
\endaligned\tag3.3$$
In the case $m_1<\cdots<m_n$, we also have
$$\aligned &[x_1^k\cdots x_n^k]\|a_{ij}x_j^{m_i}\|_{1\ls i,j\ls n}
\prod_{1\ls i<j\ls n}(x_j-x_i)\times\bg(\sum_{s=1}^nx_s\bg)^{kn-\bi n2-\sum_{i=1}^n m_i}
\\&\quad\ =(-1)^{\bi n2}\f{(kn-\bi n2-\sum_{i=1}^n m_i)!}
{\prod_{i=1}^n\prod\Sb m_i<j\ls k\\j\not\in\{m_s:\, i<s\ls n\}\endSb(j-m_i)}\det(A).
\endaligned\tag3.4$$
\endproclaim

\Remark\ 3.1. Part (ii) of Theorem 3.2 follows from the first part.
\medskip

\proclaim{Theorem 3.3} Let $R$ be a commutative ring with
identity, and let $a_{ij}\in R$ for all $i,j=1,\ldots,n$. Let $k,l_1,\cs,l_n,m_1,\cs,m_n\in\N$
with $N=kn-\sum_{i=1}^n(l_i+m_i)\gs0$.

{\rm (i) (Sun [Su3, Theorem 2.2])} There holds the identity
$$\aligned&[x_1^k\cdots x_n^k]|a_{ij}x_j^{l_i}|_{1\ls i,j\ls
n}\,|x_j^{m_i}|_{1\ls i,j\ls n}\,(x_1+\cdots+x_n)^N
\\=&[x_1^k\cdots x_n^k]|a_{ij}x_j^{m_i}|_{1\ls i,j\ls
n}\,|x_j^{l_i}|_{1\ls i,j\ls n}\,(x_1+\cdots+x_n)^N.
\endaligned\tag3.5$$

{\rm (ii)} We also have the following symmetric identities:
$$\aligned&[x_1^k\cdots x_n^k]\|a_{ij}x_j^{l_i}\|_{1\ls i,j\ls
n}\,|x_j^{m_i}|_{1\ls i,j\ls n}\,(x_1+\cdots+x_n)^N
\\=&[x_1^k\cdots x_n^k]\|a_{ij}x_j^{m_i}\|_{1\ls i,j\ls
n}\,|x_j^{l_i}|_{1\ls i,j\ls n}\,(x_1+\cdots+x_n)^N,
\endaligned\tag3.6$$
$$\aligned&[x_1^k\cdots x_n^k]|a_{ij}x_j^{l_i}|_{1\ls i,j\ls
n}\,\|x_j^{m_i}\|_{1\ls i,j\ls n}\,(x_1+\cdots+x_n)^N
\\=&[x_1^k\cdots x_n^k]|a_{ij}x_j^{m_i}|_{1\ls i,j\ls n}\,\|x_j^{l_i}\|_{1\ls i,j\ls n}\,(x_1+\cdots+x_n)^N,
\endaligned\tag3.7$$
and
$$\aligned&[x_1^k\cdots x_n^k]\|a_{ij}x_j^{l_i}\|_{1\ls i,j\ls
n}\,\|x_j^{m_i}\|_{1\ls i,j\ls n}\,(x_1+\cdots+x_n)^N
\\=&[x_1^k\cdots x_n^k]\|a_{ij}x_j^{m_i}\|_{1\ls i,j\ls n}\,\|x_j^{l_i}\|_{1\ls i,j\ls n}\,(x_1+\cdots+x_n)^N.
\endaligned\tag3.8$$
\endproclaim

\medskip
Theorem 3.3(ii) can be proved by modifying the proof of [Su3, Theorem 2.2] slightly.

\heading {4. Proof of Theorem 1.2}\endheading

\proclaim{Lemma 4.1} Let $h,k,l,m,n$ be positive integers satisfying $(1.1)$.
Let $c_1,\ldots,c_n$ be elements of a commutative ring $R$ with identity,
and let $P(x_1,\ldots,x_n,y_1,\ldots,y_n)$ denote the polynomial
$$\prod_{1\ls i<j\ls n}(c_jx_jy_j-c_ix_iy_i)(x_j^m-x_i^m)(y_j^h-y_i^h)
\times(x_1+\cdots+x_n)^K(y_1+\cdots+y_n)^L,$$
where $K$ and $L$ are given by $(1.2)$.
Then
$$\aligned&[x_1^{k-1}\cdots x_n^{k-1}y_1^{l-1}\cdots y_n^{l-1}]P(x_1,\ldots,x_n,y_1,\ldots,y_n)
\\&\qquad\qquad=\f{K!L!}{N}\prod_{1\ls i<j\ls n}(c_j-c_i),\endaligned\tag4.1$$
where
$$N=(hm)^{-\bi n2}\prod_{r=0}^{n-1}\f{(k-1-rm)!(l-1-rh)!}{(r!)^2}\in\Z^+.\tag4.2$$
\endproclaim
\Proof. In view of Theorem 3.3(i) and Theorem 3.1(ii),
$$\align &[y_1^{l-1}\cdots y_n^{l-1}]\prod_{1\ls i<j\ls n}(c_jx_jy_j-c_ix_iy_i)(y_j^h-y_i^h)\times(y_1+\cdots+y_n)^L
\\=&[y_1^{l-1}\cdots y_n^{l-1}]|(c_jx_j)^{i-1}y_j^{i-1}|_{1\ls i,j\ls n}|y_j^{(i-1)h}|_{1\ls i,j\ls n}(y_1+\cdots+y_n)^L
\\=&[y_1^{l-1}\cdots y_n^{l-1}]|(c_jx_j)^{i-1}y_j^{(i-1)h}|_{1\ls i,j\ls n}|y_j^{i-1}|_{1\ls i,j\ls n}(y_1+\cdots+y_n)^L
\\=&(-1)^{\bi n2}\f{L!}{L_0}\|(c_jx_j)^{i-1}\|_{1\ls i,j\ls n},
\endalign$$
where
$$\align L_0=&\prod_{i=1}^n\prod\Sb (i-1)h<j\ls l-1\\j/h\not\in\{s\in\Z:\,i\ls s<n\}\endSb(j-(i-1)h)
=\prod_{i=1}^n\f{(l-1-(i-1)h)!}{\prod_{0<j\ls n-i}(jh)}
\\=&\prod_{i=1}^n\f{(l-1-(i-1)h)!}{(n-i)!h^{n-i}}
=h^{-\bi n2}\prod_{r=0}^{n-1}\f{(l-1-rh)!}{r!}.
\endalign$$
Thus, with helps of Theorem 3.3(ii) and Theorem 3.2(ii), we have
$$\align&(-1)^{\bi n2}[x_1^{k-1}\cdots x_n^{k-1}y_1^{l-1}\cdots y_n^{l-1}]P(x_1,\ldots,x_n,y_1,\ldots,y_n)
\\=&[x_1^{k-1}\cdots x_n^{k-1}]\f{L!}{L_0}
\|(c_jx_j)^{i-1}\|_{1\ls i,j\ls n}\prod_{1\ls i<j\ls n}(x_j^m-x_i^m)
\times\(\sum_{s=1}^nx_s\)^K
\\=&\f{L!}{L_0}[x_1^{k-1}\cdots x_n^{k-1}]\|c_j^{i-1}x_j^{i-1}\|_{1\ls i,j\ls n}
|x_j^{(i-1)m}|_{1\ls i,j\ls n}(x_1+\cdots+x_n)^K
\\=&\f{L!}{L_0}[x_1^{k-1}\cdots x_n^{k-1}]\|c_j^{i-1}x_j^{(i-1)m}
\|_{1\ls i,j\ls n}|x_j^{i-1}|_{1\ls i,j\ls n}(x_1+\cdots+x_n)^K
\\=&\f{L!}{L_0}(-1)^{\bi n2}\f{K!}{K_0}|c_j^{i-1}|_{1\ls i,j\ls n}
=(-1)^{\bi n2}\f{K!L!}{K_0L_0}\prod_{1\ls i<j\ls n}(c_j-c_i),
\endalign$$
where
$$K_0=\prod_{i=1}^n\prod\Sb (i-1)m<j\ls k-1\\j/m\not\in\{s\in\Z:\,i\ls s<n\}\endSb(j-(i-1)m)
=m^{-\bi n2}\prod_{r=0}^{n-1}\f{(k-1-rm)!}{r!}.\tag4.3$$
Therefore (4.1) holds with $N=K_0L_0\in\Z^+$. \qed

\medskip
\noindent{\it Proof of Theorem 1.2}. Let $f(x_1,\ldots,x_n,y_1,\ldots,y_n)$ denote the polynomial

$$\align &\prod_{1\ls i<j\ls n}(P_j(x_j)-P_i(x_i))(Q_j(y_j)-Q_i(y_i))(c_jx_jy_j-c_ix_iy_i)
\\&\ \ \times (x_1+\cdots+x_n)^{K-|S|}\prod_{a\in S}(x_1+\cdots+x_n-a)
\\&\ \ \times(y_1+\cdots+y_n)^{L-|T|}\prod_{b\in T}(y_1+\cdots+y_n-b).
\endalign$$
Then
$$\deg f\ls(m+h+2)\bi n2+|K|+|L|=(k-1+l-1)n=\sum_{i=1}^n(|A_i|-1+|B_i|-1).$$
Since $\ch(F)>\max\{K,L\}$ and $\prod_{1\ls i<j\ls n}(c_j-c_i)\not=0$, in view of Lemma 4.1 we have
$$\align &[x_1^{k-1}\cdots x_n^{k-1}y_1^{l-1}\cdots y_n^{l-1}]f(x_1,\ldots,x_n,y_1,\ldots,y_n)
\\=&[x_1^{k-1}\cdots x_n^{k-1}y_1^{l-1}\cdots y_n^{l-1}]P(x_1,\ldots,x_n,y_1,\ldots,y_n)\not=0,
\endalign$$
where $P(x_1,\ldots,x_n,y_1,\ldots,y_n)$ is defined as in Lemma 4.1.
Applying the Combinatorial Nullstellensatz we find that $f(a_1,\ldots,a_n,b_1,\ldots,b_n)\not=0$ for some
$a_1\in A_1,\ldots,a_n\in A_n,b_1\in B_1,\ldots,b_n\in B_n$.
Thus (1.4) holds, and also $a_1+\cdots+a_n\not\in S$ and $b_1+\cdots+b_n\not\in T$.
We are done.  \qed

\heading{5. Proof of Theorem 1.4}\endheading

 Non-vanishing permanents are useful in combinatorics.
For example, Alon's permanent lemma [A1] states that, if
$A=(a_{ij})_{1\ls i,j\ls n}$ is a matrix over a field $F$ with $\per(A)\not=0$, and
$X_1,\ldots,X_n$ are subsets of $F$ with cardinality $2$, then for any $b_1,\ldots,b_n\in F$ there are
$x_1\in X_1,\ldots,x_n\in X_n$ such that $\sum_{j=1}^na_{ij}x_j\not=b_i$ for all $i=1,\ldots,n$.

In contrast with [Su3, Theorem 1.2(ii)], we have the following auxiliary result.
\proclaim{Theorem 5.1} Let $A_1,\ldots,A_n$ be finite subsets of a field $F$
with $|A_1|=\cdots=|A_n|=k$,
and let $P_1(x),\ldots,P_n(x)\in F[x]$ have degree at most $m\in\Z^+$
with $[x^m]P_1(x),\ldots,[x^m]P_n(x)$ distinct.
Suppose that $k-1\gs m(n-1)$ and $\ch(F)>(k-1)n-(m+1)\bi n2$.
Then the restricted sumset
$$C=\bg\{\sum_{i=1}^n a_i:\, a_i\in A_i,\ a_i\not=a_j\ \t{for}\ i\not=j,
\ \t{and}\ \|P_j(a_j)^{i-1}\|_{1\ls i,j\ls n}\not=0\bg\}\tag5.1$$
has cardinality at least $(k-1)n-(m+1)\bi n2+1>(m-1)\bi n2$.
\endproclaim
\Proof. Assume that $|C|\ls K=(k-1)n-(m+1)\bi n2$.
Clearly the polynomial
$$\align f(x_1,\ldots,x_n):=&\prod_{1\ls i<j\ls n}(x_j-x_i)\times\|P_{j}(x_j)^{i-1}\|_{1\ls i,j\ls n}
\\&\times \prod_{c\in C}(x_1+\cdots+x_n-c)\times (x_1+\cdots+x_n)^{K-|C|}
\endalign$$
has degree not exceeding $(k-1)n=\sum_{i=1}^n(|A_i|-1)$.
Since $\ch(F)$ is greater than $K$, and those $b_i=[x^m]P_i(x)$ with $1\ls i\ls n$
are distinct, with the help of Theorem 3.2(ii) we have
$$\align &[x_1^{k-1}\cdots x_n^{k-1}]f(x_1,\ldots,x_n)
\\=&[x_1^{k-1}\cdots x_n^{k-1}]\prod_{1\ls i<j\ls n}(x_j-x_i)
\times\|b_{j}^{i-1}x_j^{(i-1)m}\|_{1\ls i,j\ls n}\(\sum_{s=1}^n x_s\)^K
\\=&(-1)^{\bi n2}\f{K!}{K_0}|b_j^{i-1}|_{1\ls i,j\ls n}=(-1)^{\bi n2}\f{K!}{K_0}\prod_{1\ls i<j\ls n}(b_j-b_i)\not=0,
\endalign$$
where $K_0$ is given by (4.3).
Thus, by the Combinatorial Nullstellensatz, $f(a_1,\ldots,a_n)\not=0$ for some
$a_1\in A_1,\ldots,a_n\in A_n$. Clearly $\sum_{i=1}^na_i\in C$ if $\|P_{j}(a_j)^{i-1}\|_{1\ls i,j\ls n}\not=0$
and $a_i\not=a_j$ for all $1\ls i<j\ls n$.
So we also have $f(a_1,\ldots,a_n)=0$ by the definition of $f(x_1,\ldots,x_n)$.
The contradiction ends our proof. \qed

\proclaim{Corollary 5.1} Let $A_1,\ldots,A_n$ and $B=\{b_1,\ldots,b_n\}$
be subsets of a field with cardinality $n$. Then there is an SDR $\{a_i\}_{i=1}^n$ of $\{A_i\}_{i=1}^n$
such that the permanent $\|(a_jb_j)^{i-1}\|_{1\ls i,j\ls n}$ is nonzero.
\endproclaim
\Proof. Simply apply Theorem 5.1 with $k=n$ and $P_j(x)=b_jx$ for $j=1,\ldots,n$. \qed

\proclaim{Lemma 5.1} Let $k,m,n\in\Z^+$ with $k-1\gs m(n-1)$. Then
$$\aligned &[x_1^{k-1}\cdots x_n^{k-1}]\prod_{1\ls i<j\ls n}(x_j-x_i)^{2m-1}(x_jy_j-x_iy_i)
\times\(\sum_{s=1}^nx_s\)^{N}
\\=&(-1)^{m\bi n2}\f{(mn)!N!}{(m!)^n n!}\prod_{r=0}^{n-1}\f{(rm)!}{(k-1-rm)!}
\times\|y_j^{i-1}\|_{1\ls i,j\ls n},
\endaligned\tag5.2$$
where $N=(k-1-m(n-1))n$.
\endproclaim
\Proof. Since both sides of (5.2) are polynomials in $y_1,\ldots,y_n$, it suffices to show
that (5.2) with $y_1,\ldots,y_n$ replaced by $a_1,\ldots,a_n\in\C$ always holds.

 By Lemma 2.1 and (2.6) of [SY], we have
 $$\align &[x_1^{k-1}\cdots x_n^{k-1}]\prod_{1\ls i<j\ls n}(x_j-x_i)^{2m-1}(a_jx_j-a_ix_i)
\times\(\sum_{s=1}^nx_s\)^{N}
\\=&\f{N!}{((k-1)!)^n}(-1)^{m\bi n2}\f{m!(2m)!\cdots(nm)!}{(m!)^nn!}
\|a_j^{i-1}\|_{1\ls i,j\ls n}\prod_{0<r<n}\prod_{s=1}^{rm}(k-s)
\\=&(-1)^{m\bi n2}\f{(mn)!N!}{(m!)^n n!}\|a_j^{i-1}\|_{1\ls i,j\ls n}\prod_{r=0}^{n-1}\f{(rm)!}{(k-1-rm)!}.
 \endalign$$
 This concludes the proof. \qed

 \medskip
 \noindent{\it Proof of Theorem 1.4}. Since $c_1,\ldots,c_n$ are distinct and $|B_1|=\cdots=|B_n|=n$,
 by Corollary 5.1 there is an SDR $\{b_i\}_{i=1}^n$ of $\{B_i\}_{i=1}^n$ such that
 $\|(b_jc_j)^{i-1}\|_{1\ls i,j\ls n}\not=0$.

 Suppose that $|S|\ls N=(k-1-m(n-1))n$. We want to derive a contradiction.
 Let $f(x_1,\ldots,x_n)$ denote the polynomial
 $$\align&\prod_{1\ls i<j\ls n}\((b_jc_jx_j-b_ic_ix_i)(x_j-x_i)^{2m-1-|S_{ij}|}\prod_{c\in S_{ij}}(x_j-x_i+c)\)
 \\&\quad\qquad\times(x_1+\cdots+x_n)^{N-|S|}\prod_{a\in S}(x_1+\cdots+x_n-a).
\endalign$$
Then $$\deg f\ls 2m\bi n2+N=(k-1)n=\sum_{i=1}^n(|A_i|-1).$$
With the help of Lemma 5.1, we have
$$\align &[x_1^{k-1}\cdots x_n^{k-1}]f(x_1,\ldots,x_n)
\\=&[x_1^{k-1}\cdots x_n^{k-1}](x_1+\cdots+x_n)^N
\prod_{1\ls i<j\ls n}(b_jc_jx_j-b_ic_ix_i)(x_j-x_i)^{2m-1}
\\=&(-1)^{m\bi n2}\f{(mn)!N!}{(m!)^n n!}\prod_{r=0}^{n-1}\f{(rm)!}{(k-1-rm)!}
\times\|(b_jc_j)^{i-1}\|_{1\ls i,j\ls n}\not=0
\endalign$$
since $\ch(F)>\max\{mn,N\}$. By the Combinatorial Nullstellensatz,
there are $a_1\in A_1,\ldots,a_n\in A_n$ such that $f(a_1,\ldots,a_n)\not=0$.
On the other hand, we do have $f(a_1,\ldots,a_n)=0$,
because $a_1+\cdots+a_n\in S$ if $a_i-a_j\not\in S_{ij}$ and $a_ib_ic_i\not=a_jb_jc_j$
for all $1\ls i<j\ls n$. So we get a contradiction. \qed

\enddocument